\documentclass{amsart} 
\usepackage{amsmath}

\usepackage{mathptmx}
\usepackage{helvet}
\usepackage{courier}
\usepackage{txfonts}
\usepackage{tikz} 
\usetikzlibrary{arrows}
\usepackage{type1cm}

\usepackage{paralist}
\usepackage{epsfig} 
\usepackage{epstopdf} 
\usepackage[colorlinks=true]{hyperref}
\hypersetup{urlcolor=blue, citecolor=red}
\allowdisplaybreaks

\newtheorem{theorem}{Theorem}[section]
\newtheorem{corollary}{Corollary}

\newtheorem{lemma}[theorem]{Lemma}
\newtheorem{proposition}{Proposition}

\newtheorem{construction}{Construction}
\newtheorem{observation}{Observation}
\newtheorem*{problem}{Problem}
\theoremstyle{definition}
\newtheorem{definition}[theorem]{Definition}

\newtheorem{example}{Example}

\newcommand{\ent}{\textup{ent}}

\begin{document}
\title{Sufficient conditions for non-zero entropy and finite relations}

\author{Iztok Bani\v c, Rene Gril Rogina,  Judy Kennedy, Van Nall}


\begin{abstract}
We introduce the notions of returns, dispersions and well-aligned sets for closed relations on compact metric spaces and then we use them to obtain non-trivial sufficient conditions for such a relation to have non-zero entropy. In addition, we give a characterization of finite relations with non-zero entropy in terms of Li-Yorke and DC2-chaos.
\end{abstract}

\thanks{This work was supported by the Slovenian Research Agency under the research program P1-0285.}

\maketitle

\def\joinrel{\mkern-3mu}
\newcommand{\varproj}{\displaystyle \lim_{\multimapinv\joinrel-\joinrel-}}

\section{Introduction}\label{s0}
In topological dynamics, the study of chaotic behaviour of a dynamical system is  often based on some  properties of continuous functions.  One of the frequently studied properties of such functions in the theory of topological dynamical systems is the entropy of a continuous function $f:X\rightarrow X$ on a compact metric space $X$, which serves as a measure of the complexity of the dynamical system. This often leads to studying the entropy of the shift map $\sigma$ on the inverse limit $\varprojlim(X,f)$.  More precisely,  suppose $X$ is a compact metric space.  If $f:X \to X$ is a continuous function,  the inverse limit space generated by $f$ is 
\begin{equation*}
 \varprojlim(X,f):=\Big\{(x_{1},x_{2},x_{3},\ldots ) \in \prod_{i=1}^{\infty} X \ | \ 
\text{ for each positive integer } i,x_{i}= f(x_{i+1})\Big\},
\end{equation*}
\noindent also abbreviated as $\varprojlim f$. 
The map $f$ on $X$ induces a natural homeomorphism $\sigma$ on $\varprojlim f$,  called the shift map, defined by 
$$\sigma(x_{1},x_{2},x_3, x_4,\ldots)=(x_{2},x_{3},x_4,\dots)$$ for each $(x_{1},x_{2},x_3, x_4,\ldots)$ in  $\varprojlim f$.  It is a well-known result that the entropy of $f$ is then equal to the entropy of $\sigma$ \cite[Proposition 5.2]{BOWEN}.

To study such inverse limits $\varprojlim f$ and shift maps $\sigma:\varprojlim f\rightarrow \varprojlim f$, the study of backward orbits of points of dynamical systems $(X,f)$ is also required; note that the inverse limit $\varprojlim f$ is the space of all backward orbits in $(X,f)$.  Such backward orbits of points are actually forward orbits of points in the dynamical system $(X,f^{-1})$, if $f^{-1}$ is well-defined.  But usually, $f^{-1}$ is not a well-defined function, therefore, a more general tool is needed to study these properties. Note that for a continuous function $f:X\rightarrow X$,  the set
$$
\Gamma(f)^{-1}=\{(y,x) \in X\times X \ | \  y=f(x)\}
$$  
is a closed relation on $X$ that describes best the dynamics of $(X,f)$ in the backward direction  when $f^{-1}$ is not well-defined.  So,  generalizing topological dynamical systems $(X,f)$ to topological dynamical systems $(X,G)$ with closed relations $G$ on $X$ by making the identification $(X,f)=(X,\Gamma(f))$ is only  natural.

Recently,  many such generalizations of dynamical systems were introduced and studied (see \cite{BEK,BEGK,CP,KW,KN,LP,LYY,LWZ,M,MRT,R}, where more references may be found).  However,  there is not much known of such dynamical systems and therefore,  there are many properties of such set-valued dynamical systems that are yet to be studied.  In \cite{BEK},  the notion of topological entropy $h(f)$ of continuous functions $f:X\rightarrow X$ on compact metric spaces $X$ was generalized to the notion of topological entropy $\ent(G)$ of closed relations $G$ on compact metric spaces $X$.  In this paper, we continue our research from \cite{BEK}.  We introduce the notions of returns, dispersions and well-aligned sets for closed relations on compact metric spaces and then use them to obtain non-trivial sufficient conditions for such relations to have non-zero entropy. In addition, we give a characterization of finite relations with non-zero entropy. We also show that, unlike topological entropy for closed relations on compact metric spaces in general, in the case of finite relations, positive entropy is equivalent to the shift map on the Mahavier product being Li-Yorke chaotic as well as equivalent to DC-2 distributional chaos for the shift map, as well as equivalent to $G$ having a $(k,\varepsilon)$-return.  

We proceed as follows.  In Section \ref{s1}, basic definitions and notation that are needed later in the paper are given and presented.   In Section \ref{s2}, the topological entropy for closed relations is defined and in addition, basic results from \cite{BEK} are presented.  In Section \ref{s3}, our first main result as well as many illustrative examples and corollaries are given and proved. In our last section, Section \ref{s4}, we restrict ourselves to  finite relations on compact metric spaces. Here, our second main result, a characterization of finite relations with non-zero entropy,  is presented and proved. 

\section{Definitions and notation} \label{s1}
First, we define some properties from the continuum theory and the theory of inverse limits  that will be used later in the paper.
\begin{definition}
Let $X$ be a compact metric space. We  always use \emph{\color{blue} $\rho$} to denote the metric on $X$.
\end{definition}
\begin{definition}
Suppose $X$ is a compact metric space. If for each positive integer $n$, $f_n:X \to X$ is a continuous function, \emph{\color{blue} the inverse limit space generated by $(f_n)$} is 
$$
 \varprojlim(X,f_n):=\Big\{(x_{1},x_{2},x_{3},\ldots ) \in \prod_{i=1}^{\infty} X \ | \ 
\text{ for each positive integer } i,x_{i}= f_i(x_{i+1})\Big\}.
$$
\end{definition}

\begin{definition}
\emph{\color{blue} A continuum} is a non-empty connected compact metric space. A continuum is \emph{\color{blue} degenerate}, if it consists of only a single point. Otherwise it is \emph{\color{blue} non-degenerate}.  \emph{\color{blue} A subcontinuum} is a subspace of a continuum which itself is also a continuum.
\end{definition}

Next, we define chainable continua (using inverse limits); see \cite[Section XII]{nadler} for more details. 
\begin{definition}
A continuum $X$ is \emph{\color{blue} chainable} if there is a sequence $(f_n)$ of continuous surjections $f_n:[0,1]\rightarrow [0,1]$ such that $X$ is homeomorphic to $\varprojlim([0,1],f_n)_{n=1}^{\infty}$. 
\end{definition}
\begin{definition}
A continuum $X$ is \emph{\color{blue} decomposable}, if there are proper sub\-con\-tinua $A$ and $B$ of $X$ ($A,B\neq X$) such that $X=A\cup B$. A continuum is \emph{\color{blue} indecomposable,} if it is not decomposable. A continuum is \emph{\color{blue} hereditarily indecomposable}, if each of its subcontinua is indecomposable. 
\end{definition}
\begin{definition}
\emph{\color{blue} A pseudoarc} is any non-degenerate  hereditarily indecomposable chainable continuum.
\end{definition}
Bing showed in \cite{bing} that any two pseudoarcs are homeomorphic. 

Next, we present basic definitions and well-known results about closed relations and Mahavier products.

%
%

\begin{definition}
Let $X$ and $Y$ be metric spaces, and let $f:X\rightarrow Y$ be a function.  We use  $\Gamma(f)=\{(x,y)\in X\times Y \ | \ y=f(x)\}$
to denote \emph{ \color{blue}  the graph of the function $f$}.
\end{definition}

\begin{definition}
Let $X$ be a compact metric space and let $G\subseteq X\times X$ be a relation on $X$. If $G$ is closed in $X\times X$, then we say that $G$ is  \emph{\color{blue}  a closed relation on $X$}.  
\end{definition}

\begin{definition}
Let $X$  be a set and let $G$ be a relation on $X$.  Then we define  
$G^{-1}=\{(y,x)\in X\times X \ | \ (x,y)\in G\}$
to be \emph{\color{blue} the inverse relation of the relation $G$ on $X$}.
\end{definition}
\begin{definition}
Let $X$ be a compact metric space and let $G$ be a closed relation on $X$. Then we call
$$
\star_{i=1}^{m}G^{-1}=\Big\{(x_1,x_2,x_3,\ldots ,x_{m+1})\in \prod_{i=1}^{m+1}X \ | \ \textup{ for each } i\in \{1,2,3,\ldots ,m\}, (x_{i+1},x_i)\in G\Big\}
$$
for each positive integer $m$, \emph{\color{blue} the $m$-th Mahavier product of $G^{-1}$},  and
$$
\star_{i=1}^{\infty}G^{-1}=\Big\{(x_1,x_2,x_3,\ldots )\in \prod_{i=1}^{\infty}X \ | \ \textup{ for each positive integer } i, (x_{i+1},x_i)\in G\Big\}
$$
\emph{\color{blue} the infinite  Mahavier product of $G^{-1}$}.
\end{definition}
\begin{definition}
Let $X$ be a compact metric space and let $G$ be a closed relation on $X$.  The function  
$$
\sigma : \star_{n=1}^{\infty}G^{-1} \rightarrow \star_{n=1}^{\infty}G^{-1},
$$
 defined by 
$$
\sigma (x_1,x_2,x_3,x_4,\ldots)=(x_2,x_3,x_4,\ldots)
$$
for each $(x_1,x_2,x_3,x_4,\ldots)\in \star_{n=1}^{\infty}G^{-1}$, 
is called \emph{ \color{blue}   the shift map on $\star_{n=1}^{\infty}G^{-1}$}.  
\end{definition}

\section{Topological entropy of closed relations on compact metric spaces}\label{s2}

In this section we will summarize the generalization of topological entropy to closed relations on a compact metric space introduced in \cite{BEK}.

\begin{definition}
Let $X$ be a compact metric space and let $\mathcal S$ be a family of subsets of $X$. We use  \emph{ \color{blue}   $|\mathcal S|$} to denote the cardinality of $\mathcal S$.
\end{definition}
\begin{definition}
Let $X$ be a compact metric space and let $\mathcal S$ be a family of subsets of $X$. For each positive integer $n$, we use  $\mathcal S^n$ to denote the family
$$
\mathcal S^n=\{S_1\times S_2\times S_3\times \ldots \times S_n \ | \ S_1,S_2,S_3,\ldots,S_n\in \mathcal S\}.
$$ 
We call the elements $S_1\times S_2\times S_3\times \ldots \times S_n$ of $\mathcal S$ \emph{ \color{blue}   the $n$-boxes (generated by the family $\mathcal S$)}.
\end{definition}
\begin{definition}
Let $X$ be a compact metric space and let $\mathcal U$ be a non-empty open cover for $X$. We use \emph{\color{blue} $N(\mathcal U)$} to denote
$$
N(\mathcal U)=\min\{|\mathcal V| \ | \ \mathcal V \text{ is a non-empty finite subcover of } \mathcal U\}.
$$
\end{definition}
\begin{definition}
Let $X$ be a compact metric space, let $K$ be a closed subset of the product $\prod_{i=1}^{n}X$, and let $\mathcal U$ be a non-empty family of open subsets of $\prod_{i=1}^{n}X$ such that $K\subseteq \bigcup \mathcal U$. We use \emph{\color{blue} $N(K,\mathcal U)$} to denote
$$
N(K,\mathcal U)=\min\Big\{|\mathcal V| \ | \ \mathcal V \text{ is a non-empty subfamily of } \mathcal U \text{ such that }  K\subseteq \bigcup \mathcal V\Big\}.
$$ 
\end{definition}

\begin{theorem}\label{limit exists}\cite[Theorem 3.7]{BEK} 
Let $X$ be a compact metric space, let $G$ be a  closed relation on $X$, and let $\alpha $ be a non-empty open cover for $X$.  Then the limit
$$
\lim_{m\to \infty}\frac{\log N(\star_{i=1}^{m}G^{-1},\alpha^{m+1})}{m}
$$
exists.
\end{theorem}

\begin{definition}
Let $X$ be a compact metric space, let $G$ be a  closed relation on $X$, and let $\alpha $ be a non-empty open cover for $X$.  We define \emph{ \color{blue}   the entropy of $G$ with respect to the open cover $\alpha $} by
$$
\ent(G,\alpha )=\lim_{m\to \infty}\frac{\log N(\star_{i=1}^{m}G^{-1},\alpha^{m+1})}{m}.
$$
\end{definition}

\begin{definition}
Let $X$ be a metric space and let $\mathcal S$ and $\mathcal T$ be families of subsets of $X$. We say that the family $\mathcal S$ \emph{ \color{blue}   refines} the family $\mathcal T$, if for each $S\in \mathcal S$ there is $T\in \mathcal T$ such that $S\subseteq T$.  The notation 
$$
\mathcal T\leq \mathcal S
$$
means that the family $\mathcal S$ refines the family $\mathcal T$.
\end{definition}
\begin{proposition}\label{prop1}\cite[Proposition 1] {BEK}
Let $X$ be a compact metric space and let $G$ be a  closed relation on $X$. For all non-empty open covers $\alpha$ and $\beta$, 
$$
\alpha\leq \beta \Longrightarrow \ent(G,\alpha )\leq \ent(G,\beta ).
$$
\end{proposition}

\begin{proposition}\label{prop2} \cite[Proposition 2]{BEK}
Let $X$ be a compact metric space and let $\alpha$ be a non-empty open cover for $X$.  For all  closed relations $H$ and $G$  on $X$,
$$
H\subseteq G \Longrightarrow \ent(H,\alpha )\leq \ent(G,\alpha ).
$$
\end{proposition}

\begin{definition}
Let $X$ be a compact metric space, let $G$ be a  closed relation on $X$, and let 
$$
E=\{\ent(G,\alpha ) \ | \ \alpha \textup{ is a non-empty open cover for } X\}.
$$
  We define \emph{ \color{blue}   the entropy of $G$} by
$$
\ent(G)=\begin{cases}
				\sup (E)\text{;} & G\neq \emptyset \textup{ and } E \textup{ is bounded in } \mathbb R \\
				\infty\text{;} & G\neq \emptyset \textup{ and } E \textup{ is not bounded in } \mathbb R.
			\end{cases}
$$
\end{definition}

The following three theorems from \cite{BEK} summarize what we need to know about $\ent(G)$.

\begin{theorem}\label{thm2thm}\cite[Theorem 3.11]{BEK}
Let $X$ be a compact metric space.  For all closed relations $H$ and $G$ on  $X$, 
$$
H\subseteq G \Longrightarrow \ent(H)\leq \ent(G).
$$
\end{theorem}

\begin{theorem}\label{thm3thm}\cite[Theorem 3.12]{BEK} 
Let $X$ be a compact metric space  and let $G$ be a closed relation on $X$.  Then 
$$
\ent(G^{-1}) = \ent(G).
$$
\end{theorem}

In  \cite{BEK} it is shown that the entropy of closed relations on $X$ is a generalization of the well-known topological entropy of  continuous functions $f:X\rightarrow X$.   For a continuous function $f:X\rightarrow X$ on the compact metric space $X$ the entropy of $f$ is usually denoted $h(f)$. To suit the purposes of this paper it is enough to note that it is shown in \cite[Theorem 3.19]{BEK} that $h(f)=\ent(\Gamma(f))$. For more information on $h(f)$ see \cite{BEK} and \cite{Walters book}. 
Finally, we need Theorem \ref{the same},  also from \cite{BEK}, where the following notation is used.
\begin{definition}
For a compact metric space $X$,   we use 
$p_1:X\times X\rightarrow X$ and $p_2:X\times X\rightarrow X$
to denote \emph{\color{blue} the standard projections} defined by
$$
p_1(s,t)=s  \textup{ and } p_2(s,t)=t
$$
for all $(s,t)\in X\times X$.
\end{definition}

\begin{theorem}\label{the same}\cite[Theorem 3.18]{BEK}.
Let $X$ be a compact metric space, let $G$ be a closed relation on $X$ such that $p_1(G)\subseteq p_2(G)$ and let $\sigma$ be the shift map on $\star_{i=1}^{\infty}G^{-1}$. Then 
$$
\ent(G)=h(\sigma).
$$
\end{theorem}

\section{Returns and dispersions}\label{s3}
In this section,  we present returns and dispersions for closed relations on compact metric spaces and use them to obtain non-trivial sufficient conditions for a relation on a compact metric space to have non-zero entropy.

First, we introduce the notion of a $(k,\varepsilon)$-return on a set. 

 \begin{definition} Let $X$ be a compact metric space, let $A$ be a non-empty subset of $X$,   let $G$ be a closed relation on $X$, let $k$ be a positive integer {such that $k\geq 2$} and let $\varepsilon>0$.  We say that \emph{\color{blue} $G$ has a $(k,\varepsilon)$-return on $A$} if for each $a\in A$ there are integers $j$ and $j'$ such that $1<j'\leq j \leq k$ and points $(x_1,x_2, \ldots x_j)$ and $(y_1,y_1, \ldots ,y_j)$ in $\star_{i=1}^{j-1} G^{-1}$ such that 
 \begin{enumerate}
 \item $a=x_1=y_1$, 
 \item $\{x_j , y_j\} \subseteq A$, and 
 \item $\rho (x_{j'}, y_{j'}) >\varepsilon$.
 \end{enumerate}
 \end{definition}
 \begin{observation}
Let $X$ be a compact metric space, let $A\subseteq X$,   let $G$ be a closed relation on $X$, let $k$ be a positive integer {such that $k\geq 2$} and let $\varepsilon>0$.   Suppose that $G$ has a $(k,\varepsilon)$-return on $A$. Then for each positive integer $n>k$, $G$ has an $(n,\varepsilon)$-return on $A$.
\end{observation}
\begin{definition} Let $X$ be a compact metric space, let $G$ be a closed relation on $X$, let $k$ be a positive integer {such that $k\geq 2$} and let $\varepsilon>0$.  We say that \emph{\color{blue} $G$ has a $(k,\varepsilon)$-return} if there is a set $A\subseteq X$ such that $G$ has a $(k,\varepsilon)$-return on $A$.
 \end{definition}

Next, we give some illustrative examples.
  \begin{example} \label{ex11}
  Let $X=[0,1]$ and let $G=\{(0,0),(0,1),(1,0)\}$. Then $G$ has a $(5,\frac{1}{2})$-return on $A=\{1\}$. To see this, observe that for 
  $$
  \mathbf x=(x_1,x_2,x_3,x_4,x_5)=(1,0,1,0,1) \textup{ and } \mathbf y=(y_1,y_2,y_3,y_4,y_5)=(1,0,0,0,1),
  $$
  $\mathbf x,\mathbf y\in \star_{i=1}^{4}G^{-1}$ are such points that 
  \begin{enumerate}
 \item $1=x_1=y_1$, 
 \item for $j=5$, $\{x_j , y_j\}=\{1\}=A \subseteq A$, and 
 \item for $j'=3$, $\rho (x_{j'}, y_{j'})=1 >\frac{1}{2}$.
 \end{enumerate}
 \end{example}
   \begin{example} \label{ex13}
  Let $X=[0,1]$ and let $G=\{(0,\frac{1}{2}),(0,1),(\frac{1}{2},0),(1,0)\}$. Then $G$ has a $(5,\frac{1}{4})$-return on $A=\{1\}$. To see this, observe that for 
  $$
  \mathbf x=(x_1,x_2,x_3,x_4,x_5)=\Big(1,0,\frac{1}{2},0,1\Big) \textup{ and } \mathbf y=(y_1,y_2,y_3,y_4,y_5)=(1,0,1,0,1),
  $$
  $\mathbf x,\mathbf y\in \star_{i=1}^{4}G^{-1}$ are such points that 
  \begin{enumerate}
 \item $1=x_1=y_1$, 
 \item for $j=5$, $\{x_j , y_j\}=\{1\}=A  \subseteq A$, and 
 \item for $j'=3$, $\rho (x_{j'}, y_{j'})=\frac{1}{2} >\frac{1}{4}$.
 \end{enumerate}
 \end{example}
 \begin{example} \label{ex12}
  Let $X=[0,1]$ and let $G=\{(t,t) \ | \ t\in [0,1]\}$. Then for every positive integer $k$ and for every $\varepsilon >0$,  $G$ does not have a $(k,\varepsilon)$-return.
 \end{example}
 
The following notation is needed for the the inductive construction of a Cantor set in $\star_{i=1}^\infty G^{-1}$ that we  call a $(k,\varepsilon)$-dispersion.  
\begin{definition}
We use \emph{\color{blue} $\mathbb{N}$} to denote the set of positive integers $\{1,2,3,\ldots \}$,  and for each positive integer $j$, we use \emph{\color{blue} $\mathbb{N}_j$} to denote the set $\{1,2,3,\ldots, j\}$. 
\end{definition}
\begin{definition}
We use \emph{\color{blue} $\Sigma_2$} to denote the set
$$
\Sigma_2= \{ \mathbf s:\mathbb{N} \rightarrow \{0,1\} \}
$$
and for each positive integer $j$, we use  \emph{\color{blue} $\Sigma_2^j$} to denote the set
 $$
 \Sigma_2^j = \{ \mathbf s:\mathbb{N}_j \rightarrow \{0,1\} \}.
 $$
 \end{definition}
 {
 \begin{observation}
 The set $\Sigma_2$ can be identified with the set 
$$
\{(s_1,s_2,s_3,\ldots) \ | \ \textup{for each positive integer } k, s_k\in\{0,1\}\}
$$
 of sequences of $0$'s and $1$'s, and for each positive integer $j$,  the set   $\Sigma_2^j$ can be identified with the set
 $$
\{(s_1,s_2,s_3,\ldots,s_j) \ | \ \textup{for each } k\in \{1,2,3,\ldots,j\}, s_k\in\{0,1\}\}
 $$
 of $j$-tuples  of $0$'s and $1$'s. Therefore, when writing $\mathbf s\in \Sigma_2$, we mean
 $$
 \mathbf s=(s_1,s_2,s_3,\ldots) 
 $$
 for some sequence $(s_1,s_2,s_3,\ldots)$  of $0$'s and $1$'s, and when writing $\mathbf s\in \Sigma_2^j$, we mean
 $$
 \mathbf s=(s_1,s_2,s_3,\ldots,s_j) 
 $$
 for some $j$-tuple $(s_1,s_2,s_3,\ldots,s_j)$  of $0$'s and $1$'s.
 \end{observation}
 }
 \begin{definition}
Let $X$ be a compact metric space and let $G$ be a closed relation on $X$.  Also, let $m,n\in \mathbb N$ and let $(x_1,x_2, x_3,\ldots ,x_{n+1}) \in \star_{i=1}^n G^{-1}$ and $(y_1,y_2,y_3, \ldots , y_{m+1} )\in \star_{i=1}^m G^{-1}$ be such that  $x_{n+1}=y_1$.  Then  we define \emph{\color{blue} $(x_1, x_2, \ldots , x_{n+1})\star (y_1, y_2, \ldots , y_{m+1})$} by 
$$
(x_1, x_2, \ldots , x_{n+1})\star (y_1, y_2, \ldots , y_{m+1}) = (x_1, x_2,\ldots , x_{n+1}, y_2, \ldots , y_{m+1}).
$$
\end{definition}
{
\begin{observation}
Let $X$ be a compact metric space and let $G$ be a closed relation on $X$.  Also, let $m,n\in \mathbb N$ and let $(x_1,x_2, x_3,\ldots ,x_{n+1}) \in \star_{i=1}^n G^{-1}$ and $(y_1,y_2,y_3, \ldots , y_{m+1} )\in \star_{i=1}^m G^{-1}$ be such that  $x_{n+1}=y_1$.  Then 
$(x_1, x_2, \ldots , x_{n+1})\star (y_1, y_2, \ldots , y_{m+1}) \in \star_{i=1}^{n+m} G^{-1}$.
\end{observation}}
We also need the following definition, where a special kind of projection is defined.
\begin{definition}
Let $X$ be a compact metric space.  For each positive integer $k$,  we use $\pi_k:\prod_{i=1}^{\infty}X\rightarrow X$ to denote \emph{ \color{blue}  the $k$-th standard projection} from $\prod_{i=1}^{\infty}X$ to $X$, defined by 
$$
\pi_k(x_1,x_2,x_3,\ldots)=x_k
$$
for all $(x_1,x_2,x_3,\ldots)\in \prod_{i=1}^{\infty}X$.  
For each positive integer $n$  and for each  $k\in \{1,2,3,\ldots, n\}$, we also use use $\pi_k:\prod_{i=1}^{n}X\rightarrow X$ to denote \emph{ \color{blue}  the $k$-th standard projection} from $\prod_{i=1}^{n}X$ to $X$, defined by 
$$
\pi_k(x_1,x_2,x_3,\ldots,x_n)=x_k
$$
for all $(x_1,x_2,x_3,\ldots,x_n)\in \prod_{i=1}^{n}X$. 
\end{definition}
In Lemma \ref{con} we use notation defined in the following definition.
\begin{definition}
Let $X$ be a set, let $(i_k)$ be a sequence of positive integers, and let $\mathbf x_k=(x_{k,1},x_{k,2},x_{k,3},\ldots ,x_{k,i_k})\in X^{i_k}$ for each positive integer $k$. We define
$$
\mathbf x_1\oplus \mathbf x_2\oplus \mathbf x_3\oplus\ldots \oplus \mathbf x_n =\oplus_{k=1}^{n}\mathbf x_k=
$$
$$
(x_{1,1},x_{1,2},x_{1,3},\ldots ,x_{1,i_1},x_{2,1},x_{2,2},x_{2,3},\ldots ,x_{2,i_2},\ldots, x_{n,1},x_{n,2},x_{n,3},\ldots ,x_{n,i_n})
$$
and
$$
\mathbf x_1\oplus \mathbf x_2\oplus \mathbf x_3\oplus\ldots =\oplus_{k=1}^{\infty}\mathbf x_k=(x_{1,1},x_{1,2},x_{1,3},\ldots ,x_{1,i_1},x_{2,1},x_{2,2},x_{2,3},\ldots ,x_{2,i_2},\ldots)
$$
\end{definition}
\begin{construction}\label{con}
Let $X$ be a compact metric space, let $A\subseteq X$,   let $G$ be a closed relation on $X$, let $k$ be a positive integer  such that $k\geq 2$ and let $\varepsilon>0$. Suppose that $G$ has a $(k,\varepsilon)$-return on $A$. Then we construct a function
$$
\Psi : \bigcup_{n\in \mathbb{N}}\Sigma_2^n  \rightarrow \bigcup_{j\in \mathbb{N}_k} \star _{i=1}^{j-1} G^{-1}
$$ 
using induction on $n$ as follows.
\begin{itemize}
\item {\color{blue}$n=1$.}

\noindent  Let $\mathbf s^0=(0)\in \Sigma_2^1$ and $\mathbf s^1=(1)\in \Sigma_2^1$. Also, let $x \in A$ be any element.  Since  $G$ has a $(k,\varepsilon)$-return on $A$, there are positive integers $j$ and $j'$ such that  
$$
1<j'\leq j\leq k
$$
and  two elements $\mathbf x^0$ and $\mathbf x^1$ of $\star_{i=1}^{j-1} G^{-1}$ such that 
\begin{enumerate}
\item $\pi_1(\mathbf x^0)=\pi_1(\mathbf x^1)=x$, 
\item $\pi_{j}(\mathbf x^0), \pi_{j}(\mathbf x^1) \in  A$, and 
\item $\rho(\pi_{j'}(\mathbf x^0), \pi_{j'}(\mathbf x^1))>\varepsilon$.
\end{enumerate}
Choose and fix such integers $j$ and $j'$, and such points $\mathbf x^0$ and $\mathbf x^1$.
Then let 
$$
\Psi(\mathbf s^0)=\mathbf x^0 \textup{ and } \Psi(\mathbf s^1)=\mathbf x^1.
$$
\item {\color{blue}${n=2}$.}

\noindent For each $s\in \{0,1\}$, let $\mathbf s^{s0}=(s,0)\in \Sigma_2^2$ and $\mathbf s^{s1}=(s,1)\in \Sigma_2^2$. Also, for each $s\in \{0,1\}$, let $x_{s}=\pi_{j}(\mathbf x^{s}) \in A$.  Since  $G$ has a $(k,\varepsilon)$-return on $A$, it follows that for each $s\in \{0,1\}$ there are positive integers $j_s$ and $j_s'$ such that  
$$
1<j_s'\leq j_s\leq k
$$
and  two elements $\mathbf x^{s0}$ and $\mathbf x^{s1}$ of $\star_{i=1}^{j_s-1} G^{-1}$ such that 
\begin{enumerate}
\item $\pi_1(\mathbf x^{s0})=\pi_1(\mathbf x^{s1})=x_{s}$, 
\item $\pi_{j_s}(\mathbf x^{s0}), \pi_{j_s}(\mathbf x^{s1}) \in  A$, and 
\item $\rho(\pi_{j_s'}(\mathbf x^{s0}), \pi_{j_s'}(\mathbf x^{s1}))>\varepsilon$.
\end{enumerate}
For each $s\in \{0,1\}$, choose and fix such integers $j_s$ and $j_s'$, and such points $\mathbf x^{s0}$ and $\mathbf x^{s1}$.
Then let 
$$
\Psi(\mathbf s^{s0})=\mathbf x^{s0} \textup{ and } \Psi(\mathbf s^{s1})=\mathbf x^{s1}
$$
for  each $s\in \{0,1\}$. 
%
%
\item {\color{blue} Induction assumption.} 

\noindent Let $n$ be a positive integer and assume that for all $s_1,s_2,s_3,\ldots,s_n\in\{0,1\}$ we have already constructed the points $\mathbf x^{s_1s_2s_3\ldots s_n}\in \Sigma_2^n$, points $x_{s_1s_2s_3\ldots s_{n-1}}\in A$ and positive integers $j_{s_1s_2s_3\ldots s_{n-1}}$ and $j_{s_1s_2s_3\ldots s_{n-1}}'$ such that 
$$
1<j_{s_1s_2s_3\ldots s_{n-1}}\leq j_{s_1s_2s_3\ldots s_{n-1}}'\leq k
$$
 and that we have defined $\Psi$ by
$$
\Psi(\mathbf s^{s_1s_2s_3\ldots s_n})=\mathbf x^{s_1s_2s_3\ldots s_n}
$$
so that it is  well-defined on $\bigcup_{i \in \{1, 2, 3,...,n\}}\Sigma_2^i$.

\item {\color{blue} ${n\rightarrow n+1}$.}

\noindent For all $s_1,s_2,s_3,\ldots,s_n\in\{0,1\}$, let 
$$
\mathbf s^{s_1s_2s_3\ldots s_n0}=(s_1,s_2,s_3,\ldots,s_n,0)\in \Sigma_2^{n+1}
$$
and 
$$
\mathbf s^{s_1s_2s_3\ldots s_n1}=(s_1,s_2,s_3,\ldots,s_n,1)\in \Sigma_2^{n+1},
$$
and let 
 $$
 x_{s_1s_2s_3\ldots s_n}=\pi_{j_{s_1s_2s_3\ldots s_{n-1}}}(\mathbf x^{s_1s_2s_3\ldots s_n}) \in A.
 $$
   Since  $G$ has a $(k,\varepsilon)$-return on $A$, it follows that for all $s_1,s_2,s_3,\ldots,s_n\in \{0,1\}$ there are positive integers $j_{s_1s_2s_3\ldots s_n}$ and $j_{s_1s_2s_3\ldots s_n}'$ such that  
$$
1<j_{s_1s_2s_3\ldots s_n}'\leq j_{s_1s_2s_3\ldots s_n}\leq k
$$
and two elements $\mathbf x^{s_1s_2s_3\ldots s_n0}$ and $\mathbf x^{s_1s_2s_3\ldots s_n1}$ of $\star_{i=1}^{j_{s_1s_2s_3\ldots s_n}-1} G^{-1}$ such that 
\begin{enumerate}
\item $\pi_1(\mathbf x^{s_1s_2s_3\ldots s_n0})=\pi_1(\mathbf x^{s_1s_2s_3\ldots s_n1})=x_{s_1s_2s_3\ldots s_n}$, 
\item $\pi_{j_{s_1s_2s_3\ldots s_n}}(\mathbf x^{s_1s_2s_3\ldots s_n0}), \pi_{j_{s_1s_2s_3\ldots s_n}}(\mathbf x^{s_1s_2s_3\ldots s_n1}) \in  A$, and 
\item $\rho(\pi_{j_{s_1s_2s_3\ldots s_n}'}(\mathbf x^{s_1s_2s_3\ldots s_n0}), \pi_{j_{s_1s_2s_3\ldots s_n}'}(\mathbf x^{s_1s_2s_3\ldots s_n1}))>\varepsilon$.
\end{enumerate}
For all $s_1,s_2,s_3,\ldots ,s_n\in \{0,1\}$, choose and fix such integers $j_{s_1s_2s_3\ldots s_n}$ and $j_{s_1s_2s_3\ldots s_n}'$, and such points $\mathbf x^{s_1s_2s_3\ldots s_n0}$ and $\mathbf x^{s_1s_2s_3\ldots s_n1}$ and let 
$$
\Psi(\mathbf s^{s_1s_2s_3\ldots s_n0})=\mathbf x^{s_1s_2s_3\ldots s_n0} \textup{ and } \Psi(\mathbf s^{s_1s_2s_3\ldots s_n1})=\mathbf x^{s_1s_2s_3\ldots s_n1}
$$
for  each $s\in \{0,1\}$. 
\end{itemize}
\end{construction}
\begin{definition}
Let $X$ be a compact metric space, let $A\subseteq X$,   let $G$ be a closed relation on $X$, let $k$ be a positive integer  such that $k\geq 2$ and let $\varepsilon>0$ such that $G$ has a $(k,\varepsilon)$-return on $A$. If  $\Psi : \bigcup_{n\in \mathbb{N}}\Sigma_2^n  \rightarrow \bigcup_{j\in \mathbb{N}_k} \star _{i=1}^{j-1} G^{-1}$
is a function constructed using Construction \ref{con}, then we say that $\Psi$ is \emph{\color{blue} a $(G,A,k,\varepsilon)$-return function}.
\end{definition}
\begin{definition}
Let $X$ be a compact metric space.  For all positive integers $k$ and $\ell$ such that $k\leq \ell$,  we use $\pi_{[k,\ell]}:\prod_{i=1}^{\infty}X\rightarrow \prod_{i=k}^{\ell}X$ to denote \emph{ \color{blue}  the standard projection} that is defined by
$$
\pi_{[k,\ell]}(x_1,x_2,x_3,\ldots,x_k,x_{k+1},x_{k+2},\ldots, x_{\ell},x_{\ell +1},x_{\ell +2},\ldots)=(x_k,x_{k+1},x_{k+2},\ldots, x_{\ell}).
$$
\end{definition}
\begin{definition}
Let $X$ be a compact metric space, let $A\subseteq X$,   let $G$ be a closed relation on $X$, let $k$ be a positive integer  such that $k\geq 2$ and let $\varepsilon>0$ such that $G$ has a $(k,\varepsilon)$-return on $A$. For any $(G,A,k,\varepsilon)$-return function
$\Psi : \bigcup_{n\in \mathbb{N}}\Sigma_2^n  \rightarrow \bigcup_{j\in \mathbb{N}_k} \star _{i=1}^{j-1} G^{-1}$, 
we define the function {\color{blue}
$
B_{\Psi}: \Sigma_2  \rightarrow  \star_{i=1}^\infty G^{-1}
$ }
by 
$$
B_{\Psi}(\mathbf s)= \Psi(\pi_1(\mathbf s)) \star \Psi(\pi_{[1,2]}(\mathbf s)) \star \Psi(\pi_{[1,3]}(\mathbf s)) \star \ldots 
$$
for any  $\mathbf s\in \Sigma_2$. We also define the set \emph{\color{blue}$S_{\Psi}$} by 
$
S_{\Psi}=\{ B_{\Psi}(\mathbf s)\,\, |\,\, \mathbf s\in \Sigma_2 \}.
$ 
\end{definition}
\begin{definition} 
Let $X$ be a compact metric space, let $A\subseteq X$,   let $G$ be a closed relation on $X$, let $k$ be a positive integer  such that $k\geq 2$ and let $\varepsilon>0$ such that $G$ has a $(k,\varepsilon)$-return on $A$. \emph{\color{blue} A $(k,\varepsilon)$-dispersion for $G$ on the set $A$} is the triple 
$
(\Psi, B_{\Psi},S_{\Psi}),
$
where  $\Psi$ is a $(G,A,k,\varepsilon)$-return function.
\end{definition}
 
\begin{observation}\label{mimika}
Let $X$ be a compact metric space, let $A\subseteq X$,   let $G$ be a closed relation on $X$, let $k$ be a positive integer  such that $k\geq 2$ and let $\varepsilon>0$ such that $G$ has a $(k,\varepsilon)$-return on $A$. Then there exists a $(k,\varepsilon)$-dispersion for $G$.
\end{observation}

 \begin{theorem}\label{steven}
 Let $X$ be a compact metric space, let $A\subseteq X$,   let $G$ be a closed relation on $X$, let $k$ be a positive integer  such that $k\geq 2$ and let $\varepsilon>0$ such that $G$ has a $(k,\varepsilon)$-return on $A$. Then there is a countable subset $A'\subseteq A$ such that $G$ has a $(k,\varepsilon )$-return on $A'$.
 \end{theorem}
 
 \begin{proof} If the closed relation $G$ on the compact set $X$ has a $(k,\varepsilon)$-return on a set $A\subseteq X$, then there is a $(k,\varepsilon)$-dispersion $(\Psi, B_{\Psi},S_{\Psi})$ for $G$ on the set $A$.  Let $$ A' = \Big\{ p_1(\Psi(g)) \,\, |\,\,  g\in \bigcup_{n\in \mathbb{N}} \Sigma_2^n \Big\}.$$ 
   It follows that $A'$ is countable  and that $G$ has a $(k,\varepsilon)$-return on $A'$.
 \end{proof}
 Theorem \ref{+entropy} is our first main result of the paper. It says that for a closed relation $G$ on a compact metric space $X$, the existence of a $(k,\varepsilon)$-return on a subset of $X$ implies that the entropy of $G$ is non-zero.  

 \begin{theorem}\label{+entropy} 
 Let $X$ be a compact metric space, let $A\subseteq X$,   let $G$ be a closed relation on $X$, let $k$ be a positive integer  such that $k\geq 2$ and let $\varepsilon>0$ such that $G$ has a $(k,\varepsilon)$-return on $A$. Then
$$
\ent(G)\geq \frac{\log (2)}{k}.
$$ 
 \end{theorem}
 
  \begin{proof}  Let $(\Psi, B_{\Psi},S_{\Psi})$ be a $(k,\varepsilon)$-dispersion for $G$ on the set $A$.   Note that for each positive integer $j$,  $ \pi_{[1,j]} (S_{\Psi}) \subseteq \star_{i=1}^{j-1} G^{-1}$.
  Let $\alpha$ be an open cover for $X$ such that for any $U\in \alpha$, the diameter of $U$ is less than $\varepsilon$.  Note that for any positive integer $m$,  the number of elements of $\alpha^{mk}$ required to cover $\pi_{[1,\ldots ,mk]}(S_{\Psi})$ is the same as the number of elements of $\Sigma_2^m$ and that $|\Sigma_2^m|=2^m$.    It follows that for each positive integer $m$,  $N(\star_{i=1}^{mk-1} G^{-1},\alpha^{mk})\geq 2^m$.
    Thus 
    $$
     \ent(G) \geq \lim_{m\rightarrow \infty} \frac{\log(  N(\star_{i=1}^{mk}G^{-1},\alpha^ {mk}))}{mk}\geq \lim_{m\rightarrow \infty } \frac{ \log(2^m) }{mk}= \lim_{m\rightarrow \infty} \frac {m\log(2)}{mk}= \frac {\log(2)}{ k}.
     $$
  \end{proof}
Next, we give some examples.
\begin{example}
By Theorem \ref{+entropy}, the closed relations $G$ from Examples \ref{ex11} and \ref{ex12} have positive entropy.
\end{example}
 \begin{example} \label{ex1}
  Let $f:[0,1]\rightarrow [0,1]$ be defined by 
  $$f(x)=
\begin{cases}
2x, \text{ for } x\leq \frac{1}{2}\\
2-2x, \text{ for } x> \frac{1}{2},
\end{cases}$$
see Figure \ref{figure1}. 
\begin{figure}[h]
	\centering
		\includegraphics[width=15em]{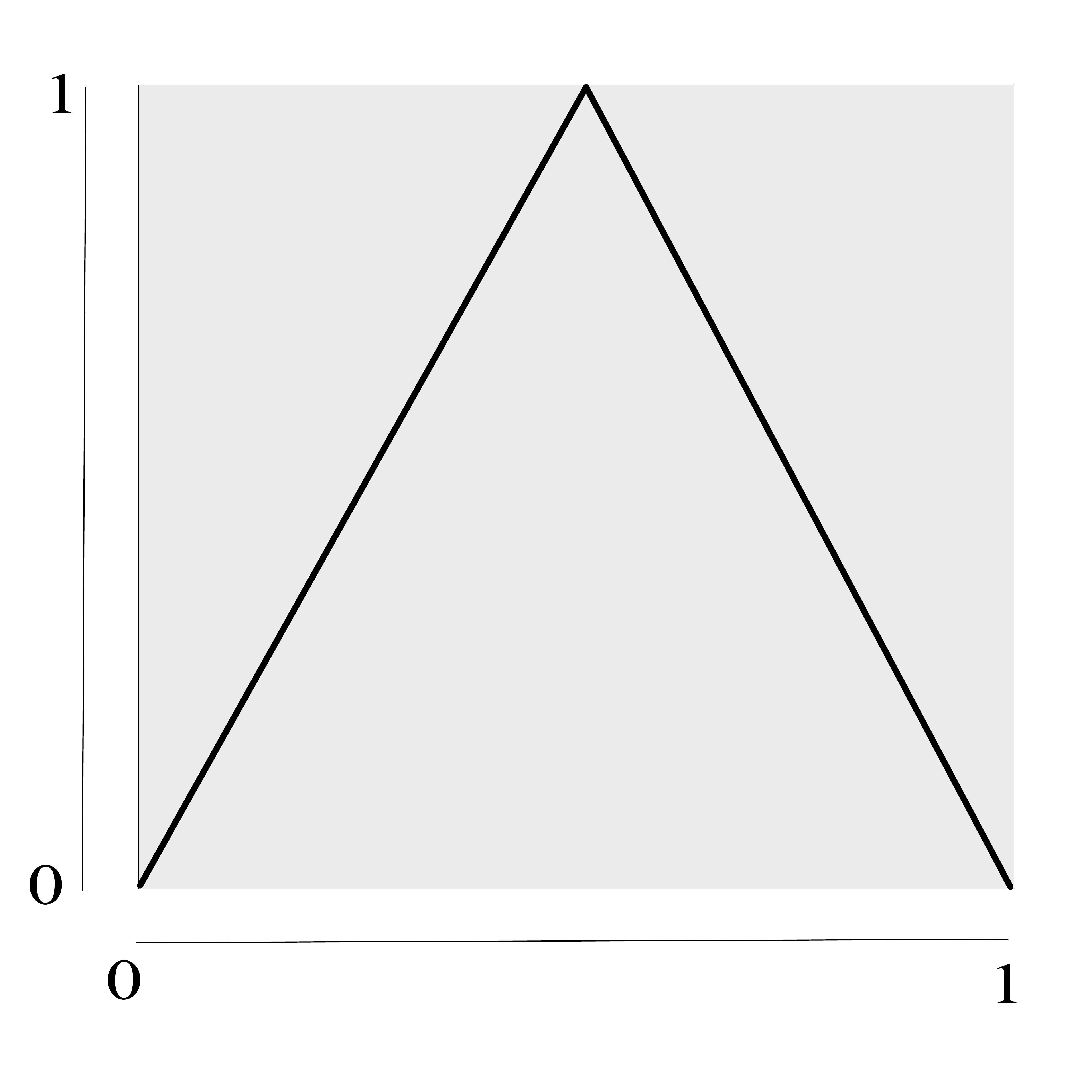}
	\caption{The graph of the function $f$ from Example \ref{ex1}}
	\label{figure1}
\end{figure} 

\noindent It is easy to check that for each $x\in [0,1]$ both $(x,\frac{1}{2}x , \frac{1}{4} x)$ and $(x, 1- \frac{1}{2}x, \frac{1}{2}+\frac{1}{4} x )$ are elements of $\star_{i=1}^2 \Gamma(f)^{-1}$. Since $|\frac{1}{4}x - ( \frac{1}{2}+\frac{1}{4} x)|=\frac{1}{2}>\frac{1}{3}$, $\Gamma(f)^{-1}$ has a $(3,\frac{1}{3})$-return on $A=[0,1]$.  It follows from Theorem \ref{+entropy} that $h(f)=\ent(\Gamma(f)^{-1})>0$.
 \end{example}
  
  \begin{example} \label{ex2}
  In \cite{BEK} it was shown that if $0<a<1$ and $0<b<1$, and 
  $$
  G=\Big\{ (x,y) \in [0,1] \times [0,1] \,\,|\,\, y=a x \text{ or } y=\frac{1}{b} x\Big\},
  $$
   then $\ent(G)\neq 0$.   We will present a much streamlined proof of this fact by showing that if $0<a\leq b<1$ then $G$ has a $(k,\varepsilon)$-return on the set $A=[ab,a]$. For the case $0 < b \leq a <1$, a similar argument shows that $G^{-1}$ has a $(k,\varepsilon)$-return on $A=[ab,b]$.
  
  So assume $0<a\leq b<1$ and $G=\{ (x,y) \in [0,1] \times [0,1] \,\,|\,\, y=a x \text{ or } y=\frac{1}{b} x \}$.  Then
  $$
  G^{-1}=\Big\{ (x,y)\in [0,1] \times [0,1] \,\,|\,\, y=\frac{1}{a} x \text{ or } y=b x \Big\},
  $$
   see Figure \ref{figure2}.
   \begin{figure}[h!]
	\centering
		\includegraphics[width=25em]{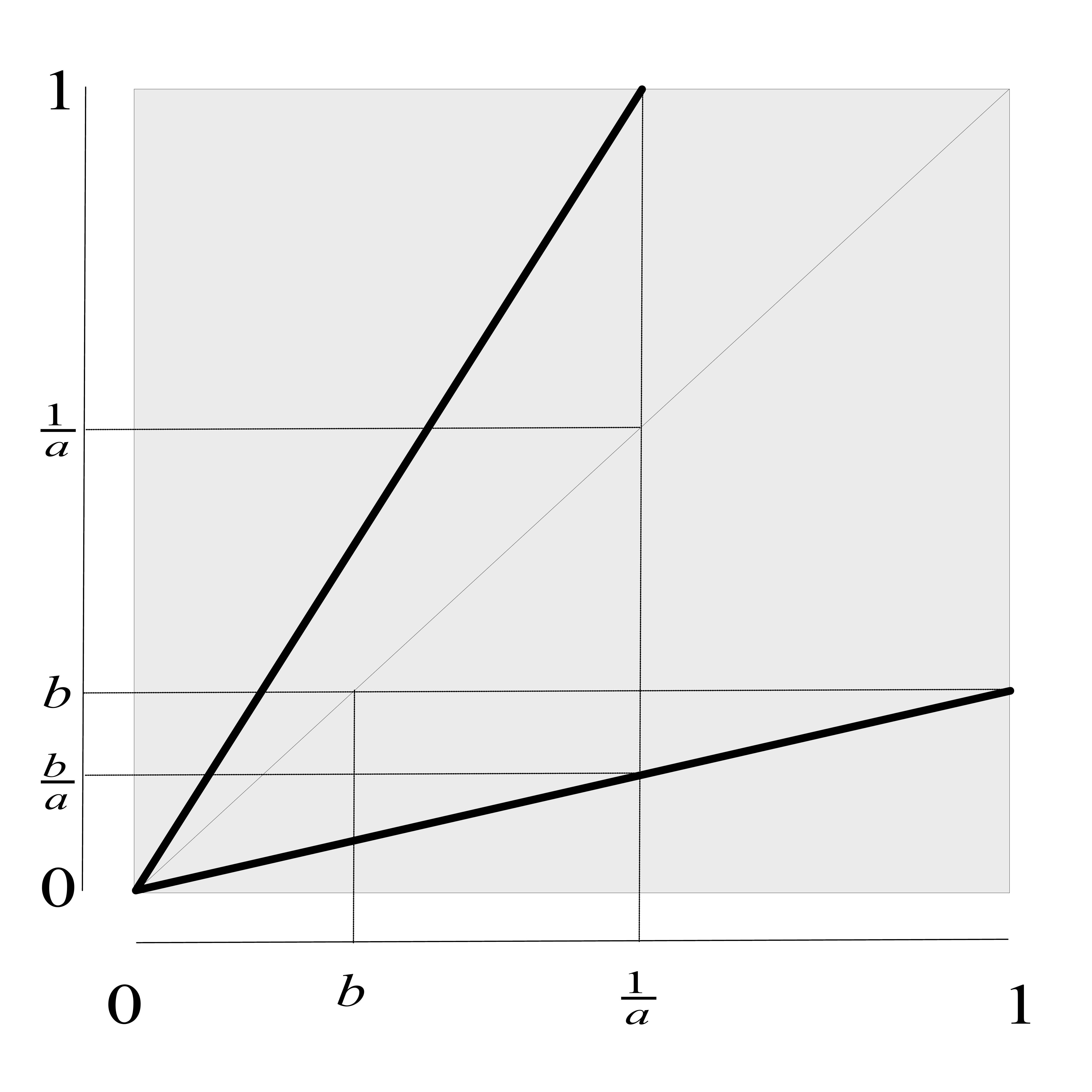}
	\caption{The relation $G^{-1}$ from Example \ref{ex2}}
	\label{figure2}
\end{figure} 
  
Note if $b>a$, then there is an $m\in \mathbb{N}$ such that $b^{m-1}>a$ and $b^m \leq a$, and thus $ab \leq b^m \leq a$.  So if $y\in (a,b]$, then there is an $m_y\in \mathbb{N}$ such that $m_y \leq m$ and $b^{m_y-1} y >a $ and $b^{m_y} y \leq a$, and thus $ab \leq b^{m_y} y\leq a$.
  
We will show that if 
$$
k=m+1  \textup{ and } \varepsilon =\frac{(\frac{1}{a}-b)ab}{2},
$$
 then $G$ has a $(k,\varepsilon)$-return on $[ab,a]$.  
  
Let $x\in [ab,a]$.  If $\frac{b}{a}x \leq a$ , then $\frac{b}{a}x \geq \frac{b}{a}ab =b\cdot b \geq ab$ and $(x,bx) \in G^{-1}$ with  $bx < x <a $ so $(bx, \frac{b}{a} x) \in G^{-1}$. It follows that 
$$
\Big(x, bx, \frac{b}{a}x\Big) \in \star_{n=1}^2 G^{-1}.
$$
  It is easy to see that also
  $$
  \Big(x,\frac{1}{a}x , \frac{b}{a} x \Big)\in \star_{n=1}^2 G^{-1}.
  $$
   Observe that $x\in [ab,a]$, $\frac{b}{a} x \in [ab,a]$ and 
$$
\Big|\frac{1}{a}x-bx\Big|=\Big(\frac{1}{a}-b\Big)x\geq \Big(\frac{1}{a}-b\Big)ab >\varepsilon.
$$
  So for $\frac{b}{a}x \leq a$ , $G$ has an $(m+1,\varepsilon)$-return on $[ab,a]$. 
 
If $\frac{b}{a}x > a$, then $\frac{b}{a} x \leq \frac{b}{a} a =b$.  So there is a positive integer $m_x$ such that 
$$
m_x\leq m \textup{ and } ab\leq b^{m_x} \frac{b}{a}x \leq a.
$$
  Now 
  $$
  \Big(x, bx, \frac{b}{a} x, b  \frac{b}{a} x, b^2  \frac{b}{a} x, \ldots , b^{m_x}  \frac{b}{a} x\Big) \in \star_{n=1}^{m_x +2} G^{-1}
  $$ 
  and 
  $$
  \Big(x, \frac{1}{a} x, \frac{b}{a} x, b  \frac{b}{a} x, b^2  \frac{b}{a} x, \ldots , b^{m_x}  \frac{b}{a} x\Big) \in \star_{n=1}^{m_x +2} G^{-1}.
  $$
   Again observe that $x\in [ab,a]$, $b^{m_x}\frac{b}{a} x \in [ab,a]$ and 
$$
\Big|\frac{1}{a}x-bx\Big|=\Big(\frac{1}{a}-b\Big)x\geq \Big(\frac{1}{a}-b\Big)ab >\varepsilon.
$$
  So for $\frac{b}{a}x > a$ , $G$ has an $(m+1,\varepsilon)$-return on $[ab,a]$. 

So, by Theorem \ref{+entropy} , it follows that $\ent(G)\neq 0$.
  \end{example}
  In the following example, we show that there are closed relations $G$ on compact metric spaces $X$ such that $\ent(G)\neq 0$ and for each non-empty $A\subseteq X$, $G$ has no $(k,\varepsilon)$-return on $A$. So,  sufficient conditions from Theorem \ref{+entropy} are not necessary conditions for non-zero entropy of $G$. Thus, we do not have a characterisation of non-zero entropy.
  \begin{example}
  Let $X$ be the pseudoarc.  Let $f:X\rightarrow X$ be a homeomorphism such that $\ent(\Gamma(f))\neq 0$.  Such a homeomorphism does exist, this was proved by J. ~Kennedy in \cite{judy1}. Since $f$ is bijective, it follows that for each non-empty $A\subseteq X$, $\Gamma(f)$ has no $(k,\varepsilon)$-return on $A$.
  \end{example}
We conclude the section by stating and proving various corollaries to Theorem \ref{+entropy}. 
  \begin{corollary} \label{juj1}
  Let $X$ be a compact metric space, let $G$ be a closed relation on $X$ and let $k$ be a positive integer.   
  If there are two  sets $J$ and $K$ in $\star_{i=1}^{k-1} G^{-1}$ with  $\rho(J,K)>0$ and such that $\pi_k(J)\cup \pi_k(K)\subseteq \pi_1(J)\cap \pi_1(K)$, then there is an $\varepsilon>0$ such that $G$ has a $(k,\varepsilon)$-return.
  \end{corollary}
  
  \begin{proof} We show that $G$ has a $(k,\varepsilon)$-return on $X$ for some $\varepsilon>0$. Let $\varepsilon=\rho (J,K)$. 
  Since $\rho (J,K)>0$, if $(x_1, x_2, \ldots ,x_k) \in J $ and $(y_1, y_2, \ldots ,y_k) \in K $ with $x_1=y_1$, then there is an integer $i\in \{2,3,4,\ldots k \}$ such that $\rho(x_i, y_i)>\varepsilon$. Let $x\in \pi_1(J)\cap \pi_1(K)$.  Then there are elements  $(x_1, x_2, \ldots ,x_k) \in J $ and $(y_1, y_2, \ldots ,y_k) \in K $ with $x=x_1=y_1$ and since $\pi_k(J)\cup \pi_k(K)\subseteq \pi_1(J)\cap \pi_1(K)$ we have $\{x_k,y_k\} \subseteq \pi_1(J)\cap \pi_1(K)$. So $G$ has a $(k,\varepsilon)$-return on $\pi_1(J)\cap \pi_1(K)$. 
  \end{proof}
   \begin{corollary} \label{juj1a}
  Let $X$ be a compact metric space, let $G$ be a closed relation on $X$ and let $k$ be a positive integer.   
  If there are two  sets $J$ and $K$ in $\star_{i=1}^{k-1} G^{-1}$ with  $\rho(J,K)>0$ and such that $\pi_k(J)\cup \pi_k(K)\subseteq \pi_1(J)\cap \pi_1(K)$, then  $\ent(G)\neq 0$.
  \end{corollary}
  \begin{proof}
  By Corollary \ref{juj1},  there is a set $A\subseteq X$, a positive integer $k$, and an $\varepsilon>0$ such that $G $ has a $(k,\varepsilon)$-return on $A$. Therefore, by Theorem \ref{+entropy},   $\ent(G)\neq 0$.
  \end{proof}

  \begin{corollary}\label{for finite} 
  Let $X$ be a compact metric space and let $G$ be a closed relation on $X$.  Also, let $k_x$ and $k_y$ be two positive integers such that $k_x,k_y>1$, let $\mathbf x\in \star_{i=1}^{k_x-1} G^{-1} $ and $\mathbf y \in \star_{i=1}^{k_y -1} G^{-1}$ be such that 
  $$
  \pi_{k_x}(\mathbf x)=\pi_1(\mathbf x)=\pi_1(\mathbf y)=\pi_{k_y}(\mathbf y),
  $$
   and let $j$ be a positive integer such that  
   $$
   1<j\leq \min\{ k_x , k_y \} \textup{ and } \pi_j(\mathbf x)\neq \pi_j(\mathbf y).
   $$
     Then there are a positive integer $k$ and an $\varepsilon>0$ such that $G $ has a $(k,\varepsilon)$-return.   
  \end{corollary}
  
  \begin{proof} 
  We show that there is a set $A\subseteq X$, a positive integer $k$, and an $\varepsilon>0$ such that $G $ has a $(k,\varepsilon)$-return on $A$. Let $\mathbf s=\mathbf x \star \mathbf y$ and $\mathbf t=\mathbf y \star \mathbf x$. Then  $\mathbf s,\mathbf t\in\star_{i=1}^{k_x+k_y-2} G^{-1} $ such that $\pi_1(\mathbf s)=\pi_1(\mathbf t)=\pi_{k_x+k_y-1}(\mathbf s)=\pi_{k_x+k_y-1}(\mathbf t)$ and $\pi_j(\mathbf s)\neq \pi_j(\mathbf t)$.  Therefore, for $A=\{ \pi_1(\mathbf s) \}$, $G$ has a $(k_x+k_y-1,\varepsilon)$-return on $A$, where $\varepsilon =\frac{1}{2} \rho(\pi_j(\mathbf s),\pi_j(\mathbf t))$.
     \end{proof}
     \begin{corollary}\label{for finite a} 
  Let $X$ be a compact metric space and let $G$ be a closed relation on $X$.  Also, let $k_x$ and $k_y$ be two positive integers such that $k_x,k_y>1$, let $\mathbf x\in \star_{i=1}^{k_x-1} G^{-1} $ and $\mathbf y \in \star_{i=1}^{k_y -1} G^{-1}$ be such that 
  $$
  \pi_{k_x}(\mathbf x)=\pi_1(\mathbf x)=\pi_1(\mathbf y)=\pi_{k_y}(\mathbf y),
  $$
   and let $j$ be a positive integer such that  
   $$
   1<j\leq \min\{ k_x , k_y \} \textup{ and } \pi_j(\mathbf x)\neq \pi_j(\mathbf y).
   $$
     Then    $\ent(G)\neq 0$.   
  \end{corollary}
  
  \begin{proof} 
  By Corollary \ref{for finite},  there is a set $A\subseteq X$, a positive integer $k$, and an $\varepsilon>0$ such that $G $ has a $(k,\varepsilon)$-return on $A$. Therefore, by Theorem \ref{+entropy},   $\ent(G)\neq 0$.
     \end{proof}
    
  In the last corollary (Corollary \ref{enkica}) to Theorem \ref{+entropy}, so-called well-aligned sets are used to detect non-zero entropy.  They form  a more visual or geometric apparatus for spotting non-zero entropy.   Before stating and proving the corollary, we give the following definitions to describe this apparatus.  
  \begin{definition}\label{dd}
Let $X$ be a compact metric space, let $G$ be a closed relation on $X$,  and let $L$ and $R$ be non-empty subsets of $G$. We say that \emph{\color{blue} the sets $L$ and $R$ are well-aligned in $G$},  if 
\begin{enumerate}
\item $p_2(L)\cap p_2(R)\neq \emptyset$,
\item there is $\varepsilon>0$ such that for all $t\in p_2(L)\cap p_2(R)$,  there are 
$$
\ell\in p_1(p_2^{-1}(t)\cap L) \text{ and } r\in p_1(p_2^{-1}(t)\cap R)
$$
 such that $\rho(r,\ell)\geq \varepsilon$,
\item $p_1(L)\cup p_1(R)\subseteq p_2(L\cup R)$, 
\item there is a positive integer $N$ such that for each $t\in p_2(L\cup R)$, there are a positive integer $i_0\leq N$ and a point 
$$
(a_1,a_2,a_3,\ldots ,a_{i_0},a_{i_0+1})\in \star_{i=1}^{i_0}G^{-1}
$$
 such that
\begin{enumerate}
\item $a_1=t$ and
\item $a_{i_0+1}\in p_2(L)\cap p_2(R)$;
\end{enumerate}
\end{enumerate}
see Figure \ref{figure3}.
 \begin{figure}[h!]
	\centering
		\includegraphics[width=32em]{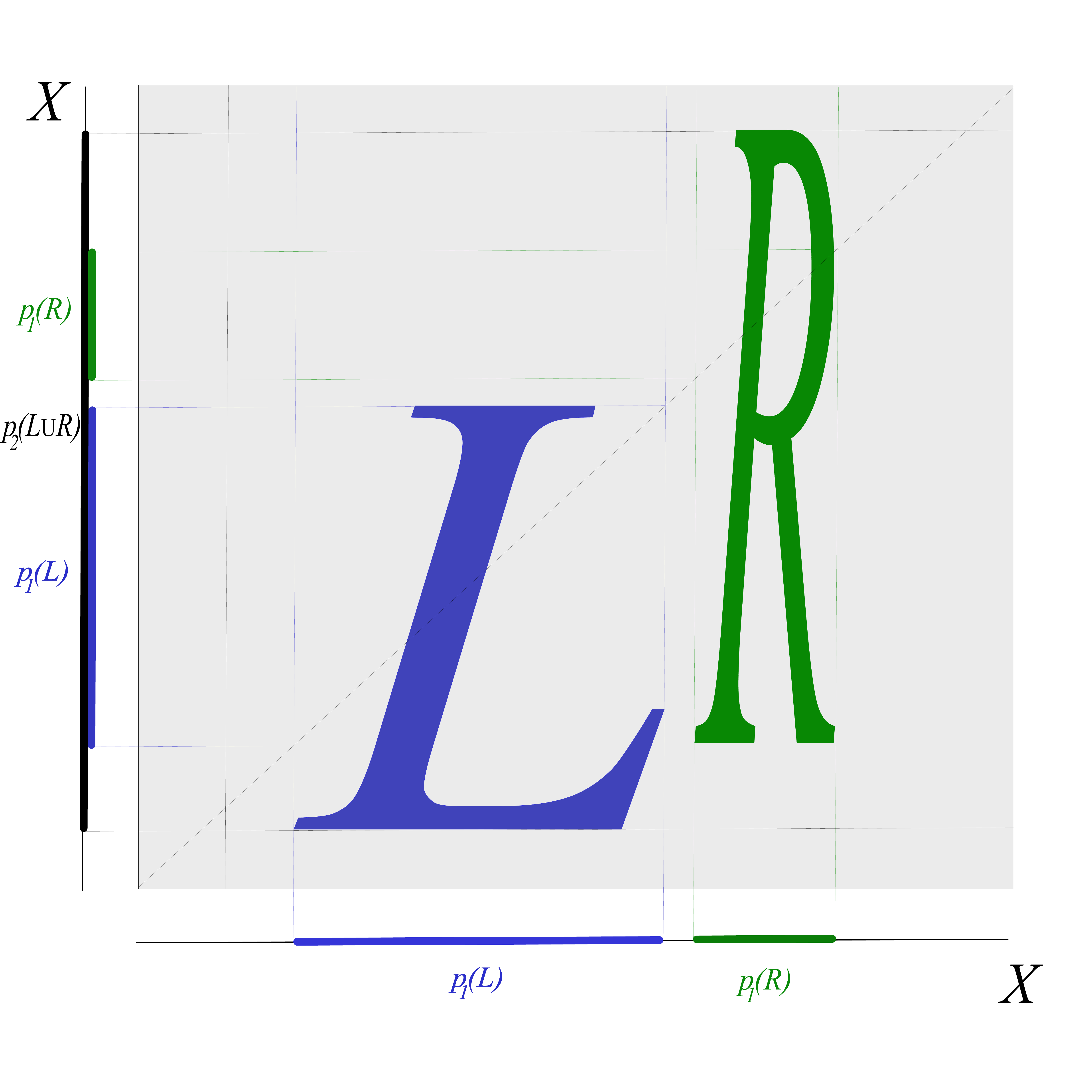}
	\caption{The sets $L$ and $R$ from Definition \ref{dd}}
	\label{figure3}
\end{figure} 
\end{definition}
\begin{example}
Let $G=\{(0,1),(0,\frac{3}{4}),(\frac{3}{4},0),(1,0)\}$ be a closed relation on $[0,1]$.
It is easy to see that $\ent(G)\neq 0$. Let $L=\{(\frac{3}{4},0)\}$, $R=\{(1,0),(0,1),(0,\frac{3}{4})\}$.  Then the sets $L$ and $R$ are well-aligned in $G$.
\end{example}
\begin{definition}
Let $X$ be a compact metric space, let $G$ be a closed relation on $X$. We say that \emph{\color{blue} the relation $G$ is well-aligned},  if there are $L,R\subseteq G$ such that the sets $L$ and $R$ are well-aligned in $G$. 
\end{definition}

\begin{corollary}\label{enkica}
Let $X$ be a compact metric space and let $G$ be a closed relation on $X$.  If $G$ or $G^{-1}$ is a well-aligned relation,  then $\ent(G)\neq 0$.
\end{corollary}

\begin{proof}
Let $L$ and $R$ be non-empty closed subsets of $G$ such that $L$ and $R$ are well-aligned in $G$; i.e.,
\begin{enumerate}
\item $p_2(L)\cap p_2(R)\neq \emptyset$,
\item there is $\varepsilon>0$ such that for all $t\in p_2(L)\cap p_2(R)$,  there are 
$$
\ell\in p_1(p_2^{-1}(t)\cap L) \text{ and } r\in p_1(p_2^{-1}(t)\cap R)
$$
 such that $\rho(r,\ell)\geq \varepsilon$,
\item $p_1(L)\cup p_1(R)\subseteq p_2(L\cup R)$, 
\item there is a positive integer $N$ such that for each $t\in p_2(L\cup R)$, there are a positive integer $i_0\leq N$ and a point 
$$
(a_1,a_2,a_3,\ldots ,a_{i_0},a_{i_0+1})\in \star_{i=1}^{i_0}G^{-1}
$$
 such that
\begin{enumerate}
\item $a_1=t$ and
\item $a_{i_0+1}\in p_2(L)\cap p_2(R)$;
\end{enumerate}
\end{enumerate}
Choose and fix such a positive integer $N$ and $\varepsilon>0$. Let $A=p_2(L\cup R)$. Then $G$ has a $(N,\frac{\varepsilon}{2})$-return on $A$. Thus, by Theorem \ref{+entropy}, $\ent(G)\neq 0$. 
\end{proof}
  
 \section{Finite relations}\label{s4}
 
Finite relations on compact metric spaces $X$ can have positive entropy. This may happen even in the case where $X$ is not finite.  We show in this final section that,  unlike topological entropy for closed relations on compact metric spaces in general, in the case of finite relations, positive entropy is equivalent to the shift map on the Mahavier product being Li-Yorke chaotic as well as equivalent to DC-2 distributional chaos for the shift map, as well as equivalent to $G$ having a $(k,\varepsilon)$-return.  
Before stating and proving our theorems, we present the following definitions.
\begin{definition}
Let $X$ be a compact metric space, let $f:X\rightarrow X$ be a continuous function, and let $x,y\in X$ such that $x\neq y$. The set $\{x,y\}$ is called \emph{\color{blue} a Li-Yorke pair for $f$} if 
$$
\liminf(\rho(f^n(x),f^n(y)))=0 \textup{ and } \limsup(\rho(f^n(x),f^n(y)))>0.
$$ 
\end{definition}

\begin{definition}
Let $X$ be a compact metric space, let $f:X\rightarrow X$ be a continuous function, and let $S\subseteq X$. We say that the set $S$ is \emph{\color{blue} a scrambled set or a Li-Yorke set in $(X,f)$}, if for all $x,y\in S$,
$$
x\neq y \Longrightarrow \{x,y\} \textup{ is a Li-Yorke pair for } f.
$$
\end{definition}
\begin{definition}
Let $X$ be a compact metric space, let $f:X\rightarrow X$ be a continuous function. The dynamical system $(X,f)$ is called \emph{\color{blue} Li–Yorke chaotic} if $X$ contains an uncountable scrambled set.
\end{definition}
The following is a well-known result.
\begin{theorem}\label{li}
Let $X$ be a compact metric space, let $f:X\rightarrow X$ be a continuous function. If $h(f)>0$, then the dynamical system $(X,f)$ is Li–Yorke chaotic. 
\end{theorem}
\begin{proof}
See the proof of \cite[2. from Corollary 2.4, page 10]{BGKM}.
\end{proof}
  See \cite{BGKM} for more references and information about Li-Yorke chaotic topological dynamical systems.
  The following lemma will be used in the proof of Theorem \ref{ch}, which is our second main result.
   \begin{lemma}    \label{for finite1} 
     Let $X$ be a compact metric space and let $G$ be a non-empty  closed relation on $X$.  If  for each positive integer $k\geq 2$ and for each $\varepsilon>0$,  $G$ has no $(k,\varepsilon)$-returns, then for each $a\in p_1(G)$, there is at most one  point 
     $$
     \mathbf x\in \star_{i=1}^{\infty}G^{-1}
     $$
 such that    
     \begin{enumerate}
     \item $\pi_1(\mathbf x)=a$ and 
     \item for infinitely many integers $n$, $\pi_n(\mathbf x)=a$.
     \end{enumerate}
     \end{lemma}
     \begin{proof}
    Let $a\in p_1(G)$ and suppose that there are points $\mathbf x_1, \mathbf x_2\in \star_{i=1}^{\infty}G^{-1}$ such that
    \begin{enumerate}
   \item  $\mathbf x_1\neq \mathbf x_2$,
   \item $\pi_1(\mathbf x_1)=\pi_1(\mathbf x_2)=a$,  
   \item for infinitely many integers $n$, $\pi_n(\mathbf x_1)=a$ and 
    \item for infinitely many integers $n$, $\pi_n(\mathbf x_2)=a$.
\end{enumerate}     
Then there is a positive integer $j$ such that 
$$
 \pi_j(\mathbf x_1)\neq \pi_j(\mathbf x_2).
 $$
 Fix such a positive integer $j$.  Also, let $k_x$ and $k_y$  be two positive integers  such that $k_x,k_y>1$, $ 1<j\leq \min\{ k_x , k_y \}$  and 
$$
\pi_{k_x}(\mathbf x_1)=\pi_{k_y}(\mathbf x_2)=a.
$$
Let 
$$
\mathbf x=\pi_{[1,k_x]}(\mathbf x_1)\in \star_{i=1}^{k_x-1} G^{-1}  \textup{ and } \mathbf y =\pi_{[1,k_y]}(\mathbf x_2)\in \star_{i=1}^{k_y -1} G^{-1}. 
$$
 Then 
  $$
  \pi_{k_x}(\mathbf x)=\pi_1(\mathbf x)=\pi_1(\mathbf y)=\pi_{k_y}(\mathbf y),
  $$
   and  $j$ is a positive integer such that  
   $$
   1<j\leq \min\{ k_x , k_y \} \textup{ and } \pi_j(\mathbf x)\neq \pi_j(\mathbf y).
   $$
     It follows from Corollary \ref{for finite} that  there do exist a positive  integer $k\geq 2$ and an $\varepsilon>0$ such that $G$ has a $(k,\varepsilon)$ return.    
     \end{proof}
        \begin{lemma}    \label{for finite2} 
     Let $X$ be a compact metric space and let $G$ be closed relation on $X$ such that $p_1(G) \subseteq p_2(G)$.  If $\ent(G)\neq 0$, then $\star_{i=1}^\infty G^{-1}$ is uncountable.
     \end{lemma}
     \begin{proof}
      Assume that $\ent(G)\neq 0$.  It follows that $G\neq \emptyset$. By Theorem \ref{the same},  the topological entropy of the shift map $\sigma$ on $\star_{i=1}^\infty G^{-1}$ is positive.  Therefore, by Theorem \ref{li}, the dynamical system $(\star_{i=1}^\infty G^{-1},\sigma)$  is Li-Yorke chaotic, which implies that there is an uncountable scrambled set in $\star_{i=1}^\infty G^{-1}$. Therefore, $\star_{i=1}^\infty G^{-1}$ is uncountable. 
     \end{proof}
  \begin{theorem} \label{ch}
Let $X$ be a compact metric space and let $G$ be a finite subset of $X \times X$ such that $p_1(G) \subseteq p_2(G)$.  The following statements are equivalent.
\begin{enumerate}
\item \label{1}   $\ent(G)\neq 0$. 
\item \label{2} There are a set $A\subseteq X$, a positive integer $k\geq 2$, and an $\varepsilon>0$ such that $G$ has a $(k,\varepsilon)$-return on $A$.
\item \label{3} There are 
\begin{enumerate}
\item positive integers $k_x$ and $k_y$ such that $k_x,k_y>1$,  
\item points $\mathbf x\in \star_{i=1}^{k_x-1} G^{-1} $ and $\mathbf y \in \star_{i=1}^{k_y -1} G^{-1}$ such that 
  $$
  \pi_{k_x}(\mathbf x)=\pi_1(\mathbf x)=\pi_1(\mathbf y)=\pi_{k_y}(\mathbf y),
  $$
   and 
   \item a positive integer $j$ such that  
   $$
   1<j\leq \min\{ k_x , k_y \} \textup{ and } \pi_j(\mathbf x)\neq \pi_j(\mathbf y).
   $$
\end{enumerate}
\item \label{4} $\star_{i=1}^\infty G^{-1}$ is uncountable.
\end{enumerate}
  \end{theorem}
  \begin{proof} 
  The implication from \ref{1} to \ref{4} is Lemma \ref{for finite2}, the implication from \ref{2} to \ref{1} is Theorem \ref{+entropy} and the implication from \ref{3} to \ref{1} is Corollary \ref{for finite a}. 
  
Now we prove the implication from \ref{4} to \ref{1}.  Assume $\ent(G)=0$.  By Theorem \ref{+entropy},  there do not exist a non-empty set $A\subseteq X$,  a positive integer $k\geq 2$, and an $\varepsilon>0$ such that $G$ has a $(k,\varepsilon)$ return on $A$.  According to Lemma \ref{for finite1}, for each $a\in p_1(G) $  there is at most one point in $\star_{i=1}^\infty G^{-1}$ with first coordinate  $a$ and in which $a$ occurs as a coordinate infinitely many times.  Let $B$ be the set of all points in $\star_{i=1}^\infty G^{-1}$ whose first coordinate occurs infinitely many times as a coordinate.  Since $G$ is finite,  $B$ is finite.  Also, for each element  $\mathbf x\in\star_{i=1}^\infty G^{-1}$, some coordinate must be repeated infinitely many times,  so there is a non-negative integer $k$ such that $\sigma^k(\mathbf x)\in B$. Since $\sigma^{-1}(\mathbf x)$ is finite for each $\mathbf x\in \star_{i=1}^\infty G^{-1}$ and since 
$$
\star_{i=1}^{\infty} G^{-1} = \bigcup _{k=0}^\infty \sigma^{-k}(B),
$$
 it follows that $\star_{i=1}^{\infty} G^{-1}$ is countably infinite or finite.
  
We prove the implication from \ref{1} to \ref{2} similarly as the implication from \ref{4} to \ref{1}:  if we assume that there does not exist a set $A\subseteq X$,  integer $k$, and $\varepsilon>0$ such that $G$ has a $(k,\varepsilon)$-return on $A$, then again Lemma \ref{for finite1} can be used as above in the proof of the implication from \ref{4} to \ref{1} to show that $\star_{i=1}^{\infty} G^{-1}$ is countable, and thus, by Lemma \ref{for finite2}, it follows that $\ent(G)=0$. 

Finally, we prove the implication from \ref{2} to \ref{3}.  Let $A\subseteq X$,  let $k\geq 2$ be a positive integer, and let $\varepsilon>0$ such that $G$ has a $(k,\varepsilon)$-return on $A$.  By Observation \ref{mimika},  there is a $(k,\varepsilon)$-dispersion $(\Psi, B_{\Psi},S_{\Psi})$ for $G$.  Fix such a dispersion.  
Since $G$ is finite,  there are at least two different sequences $\mathbf s,  \mathbf t\in \Sigma_2$ of $0$'s and $1$'s for which there are strictly increasing sequences $(i_n)$ and $(j_n)$ of positive integers such that
\begin{enumerate}
\item $i_1=1$ and $j_1=1$,
\item for each positive integer $\ell$, 
$$
\pi_{i_{\ell}}(B_{\Psi}(\mathbf s)),\pi_{j_{\ell}}(B_{\Psi}(\mathbf t))\in A 
$$
and
\item for all positive integers $\ell_1$ and $\ell_2$, 
$$
\pi_{i_{\ell_1}}(B_{\Psi}(\mathbf s))=\pi_{i_{\ell_2}}(B_{\Psi}(\mathbf s))=\pi_{j_{\ell_1}}(B_{\Psi}(\mathbf t))=\pi_{j_{\ell_2}}(B_{\Psi}(\mathbf t)). 
$$
\end{enumerate}
Since $\mathbf s\neq \mathbf t$, it follows that there is a positive integer $m$ such that 
$$
\pi_{m}(B_{\Psi}(\mathbf s))\neq \pi_{m}(B_{\Psi}(\mathbf t)).
$$
Fix such a positive integer $m$.  Then, let 
\begin{enumerate}
\item $\ell_1$ and $\ell_2$ be such positive integers that $i_{\ell_1}>m$ and $j_{\ell_2}>m$ and let 
$$
k_x=i_{\ell_1} \textup{ and } k_y=j_{\ell_2}.
$$ 
\item $\mathbf x=\pi_{[1,k_x]}(B_{\Psi}(\mathbf s))$ and $\mathbf y = \pi_{[1,k_y]}(B_{\Psi}(\mathbf t))$, and
\item $j=m$.
\end{enumerate}
Note that for constructed $k_x$, $k_y$, $\mathbf x$, $\mathbf y$ and $j$,  \ref{3} follows. This completes the proof.
  \end{proof}
 
  \begin{corollary}
 Let $X$ be a compact metric space and let $G$ be a finite subset of $X \times X$ such that $p_1(G) \subseteq p_2(G)$.  If $\ent(G)\neq 0$, then  there are elements $\mathbf s_1$ and $\mathbf s_2$ of $G$ such that 
 $$
 p_1(\mathbf s_1)=p_1(\mathbf s_2) \textup{ and } p_2(\mathbf s_1) \neq p_2 (\mathbf s_2),
 $$
  and there are elements $\mathbf t_1$ and $\mathbf t_2$ of $G$ such that  
  $$
  p_2(\mathbf t_1)=p_2(\mathbf t_2) \textup{ and } p_1(\mathbf t_1)\neq p_1(\mathbf t_2).
  $$ 
  \end{corollary}
  \begin{proof}
  Let $\ent(G)\neq 0$.  By Theorem \ref{ch},  there are
   \begin{enumerate}
\item positive integers $k_x$ and $k_y$ such that $k_x,k_y>1$,  
\item points $\mathbf x\in \star_{i=1}^{k_x-1} G^{-1} $ and $\mathbf y \in \star_{i=1}^{k_y -1} G^{-1}$ such that 
  $$
  \pi_{k_x}(\mathbf x)=\pi_1(\mathbf x)=\pi_1(\mathbf y)=\pi_{k_y}(\mathbf y),
  $$
   and 
   \item a positive integer $j$ such that  
   $$
   1<j\leq \min\{ k_x , k_y \} \textup{ and } \pi_j(\mathbf x)\neq \pi_j(\mathbf y).
   $$
\end{enumerate}
Fix such $k_x$, $k_y$, $\mathbf x$ and $\mathbf y$ and let $j_1$ be the smallest among all positive integers $j\in\{2,3,4,\ldots, k_x\}$ such that $\pi_j(\mathbf x)\neq \pi_j(\mathbf y)$.  Let 
$$
\mathbf s_1=(\pi_{j_1-1}(\mathbf x),\pi_{j_1}(\mathbf x)) \textup{ and } \mathbf s_2=(\pi_{j_1-1}(\mathbf y),\pi_{j_1}(\mathbf y)). 
 $$ 
Then $p_1(\mathbf s_1)=p_1(\mathbf s_2)$  and $p_2(\mathbf s_1) \neq p_2 (\mathbf s_2)$. This proves the first part of the claim.  Next, we prove the second part of the claim. 
Suppose that there are no elements $\mathbf t_1$ and $\mathbf t_2$ in $G$ such that  
  $$
  p_2(\mathbf t_1)=p_2(\mathbf t_2) \textup{ and } p_1(\mathbf t_1)\neq p_1(\mathbf t_2).
  $$ 
Then $G^{-1}$ is a graph of a single-valued function $f:p_2(G)\rightarrow p_1(G)$.    Note that $\varphi:p_2(G)\rightarrow \star_{i=1}^{\infty}G^{-1}$, defined by 
$$
\varphi(x)=(x,f(x),f^2(x),f^3(x),\ldots)
$$
for any $x\in p_2(G)$,  is a homeomorphism.    Since $G$ is finite, it follows that $p_2(G)$ is finite.  Therefore,   $\star_{i=1}^{\infty}G^{-1}$ is finite -- a contradiction since by Theorem \ref{ch},  $\star_{i=1}^{\infty}G^{-1}$ is uncountable. 
  This completes the proof.
  \end{proof}
  
  \begin{definition}
Let $X$ be a compact metric space, let $f:X\rightarrow X$ be a continuous function, and let $x,y\in X$ such that $x\neq y$. The set $\{x,y\}$ is called \emph{\color{blue} a DC2-pair for $f$} if 
$$
\liminf\Big(\frac{1}{n}\cdot \sum_{i=1}^{n} \rho(f^i(x),f^i(y))\Big)=0 \textup{ and } \limsup\Big(\frac{1}{n}\cdot \sum_{i=1}^{n} \rho(f^i(x),f^i(y))\Big)>0.
$$ 
\end{definition}

\begin{definition}
Let $X$ be a compact metric space, let $f:X\rightarrow X$ be a continuous function, and let $S\subseteq X$. We say that the set $S$ is \emph{\color{blue} a DC2-scrambled set in $(X,f)$}, if for all $x,y\in S$,
$$
x\neq y \Longrightarrow \{x,y\} \textup{ is a DC2-pair for }  f.
$$
\end{definition}
\begin{definition}
Let $X$ be a compact metric space, let $f:X\rightarrow X$ be a continuous function. The dynamical system $(X,f)$ is called \emph{\color{blue} DC2-chaotic} if $X$ contains an uncountable DC2-scrambled set.
\end{definition}
  
 \begin{observation}\label{iztokec}
Let $X$ be a  compact metric space and let $f:X\rightarrow X$ be a continuous function.  
 If $(X,f)$ is DC2-chaotic, then $(X,f)$ is Li-Yorke chaotic. See \cite{jana} for more information. 
  \end{observation}
Putting these facts together we see that the case where $G$ is a finite subset of $X\times X$  produces a dynamical system $(\star_{i=1}^\infty G^{-1},\sigma)$ where several forms of chaos are equivalent.  See the following corollary.  Note that it is a  known fact that in general,  Li-Yorke chaos does not imply positive entropy, see  \cite{S1} for more information.
 \begin{corollary} Suppose $X$ is a compact metric space and $G$ is a finite subset of $X \times X$ such that $p_1(G) \subseteq p_2(G)$, then the following are equivalent.
 \begin{enumerate}
 \item \label{11} $\ent(G)\neq 0$.
  \item \label{22}  $(\star_{i=1}^\infty G^{-1},\sigma)$ is Li-Yorke chaotic.
 \item \label{33} $(\star_{i=1}^\infty G^{-1},\sigma)$ has a DC2-scrambled Cantor set. 
 \end{enumerate}
  \end{corollary} 
  \begin{proof}
  First, we show the implication from \ref{11} to \ref{33}.  Suppose that $\ent(G)\neq 0$. By Theorem \ref{the same}, $h(\sigma)\neq 0$, where $\sigma$ is the shift  map on $\star_{i=1}^{\infty}G^{-1}$.  By  \cite[Theorem 1.1, p 138]{DO} $(\star_{i=1}^\infty G^{-1},\sigma)$ is DC2-chaotic and it follows from  \cite[Remark 3, p 148]{DO} that $(\star_{i=1}^\infty G,\sigma)$ has a DC2-scrambled Cantor set.   
  
  The implication from \ref{33} to \ref{22} follows from Observation \ref{iztokec}.  
  
Finally, we prove the implication from \ref{22} to \ref{11}.  Suppose that $(\star_{i=1}^\infty G^{-1},\sigma)$ is Li-Yorke chaotic. Then $\star_{i=1}^\infty G^{-1}$ contains an uncountable scrambled set. It follows that $\star_{i=1}^\infty G^{-1}$ is uncountable.  By Theorem \ref{ch},  $\ent(G)\neq 0$. 
  \end{proof}
 
 The following example, where a countable closed subset $G$ of $[0,1]\times [0,1]$,  such that 
\begin{enumerate}
\item  $\ent(G)=0$ and 
\item $\star_{i=1}^\infty G^{-1}$ is uncountable,
\end{enumerate} 
 is presented,  together with two problems, is a good place to finish the paper.
  
\begin{example} 
   Let 
   $$
   G^{-1}=\{(0,0)\}\cup \Big\{ \Big(\frac{1}{2^i},\frac{1}{2^{i+1}}\Big) \,\,|\,\,i\in \{0,1,2,3,\ldots \} \Big\}\cup 
   $$
   $$
   \Big\{ \Big(\frac{1}{2^i},\frac{1}{2^{i+1}}\Big) \,\,|\,\,i\in \{0,1,2,3,\ldots \}\Big\} \cup 
   \Big\{ \Big(\frac{1}{2^i},\frac{1}{2^{i+2}}\Big) \,\,|\,\,i\in \{0,1,2,3,\ldots \} \Big\} 
   $$
    Then $\ent(G)=0$ since $G^{-1} \subseteq\{ (x,y) \in [0,1]\times [0,1] \,\,|\,\, y\leq x \}$ (see \cite{EK}), and $\star_{i=1}^\infty G^{-1}$ is uncountable since each coordinate $a$ of an element of $\star_{i=1}^\infty G^{-1}$ can be followed by  either $\frac{1}{2} a$ or $\frac{1}{4} a$. 
\end{example}
 
 \begin{problem}
 Let $X$ be a compact metric space and let $G$ be a countable closed relation on $X$ such that $\ent(G)\neq 0$.  Is it true that either $G$ or $G^{-1}$ has a $(k,\varepsilon)$-return for some positive integer $k$ and some $\varepsilon>0$?
 \end{problem}
 \begin{problem}
 Let $X=[0,1]$ and let $f:X\rightarrow X$ be a continuous function such that $h(f)\neq 0$.  Is it true that $\Gamma(f)$ has a $(k,\varepsilon)$-return for some positive integer $k$ and some $\varepsilon>0$?
 \end{problem}

\-

\noindent I. Bani\v c\\
              (1) Faculty of Natural Sciences and Mathematics, University of Maribor, Koro\v{s}ka 160, SI-2000 Maribor,
   Slovenia; \\(2) Institute of Mathematics, Physics and Mechanics, Jadranska 19, SI-1000 Ljubljana, 
   Slovenia; \\(3) Andrej Maru\v si\v c Institute, University of Primorska, Muzejski trg 2, SI-6000 Koper,
   Slovenia\\
             {iztok.banic@um.si}           

                 \-
				
		\noindent R.  Gril Rogina\\
             Faculty of Natural Sciences and Mathematics, University of Maribor, Koro\v{s}ka 160, SI-2000 Maribor, Slovenia\\
{{rene.gril@student.um.si}       }    

                 	\-
					
  \noindent J.  Kennedy\\
             Department of Mathematics,  Lamar University, 200 Lucas Building, P.O. Box 10047, Beaumont, TX 77710 USA\\
{{kennedy9905@gmail.com}       }    
                 
           			\-
				
		\noindent V.  Nall\\
             Department of Mathematics,  University of Richmond, Richmond, VA 23173 USA\\
{{vnall@richmond.edu}       }    

   \end{document}